\documentclass[]{amsart}
\usepackage[foot]{amsaddr}
\usepackage[a4paper,top=2.5cm,bottom=3cm,inner=3cm,outer=3cm]{geometry}

\usepackage[utf8]{inputenc}
\usepackage{tikz}

\usepackage{amsmath,amsthm}
\usepackage{amsfonts,amssymb,mathrsfs}

\usepackage[bookmarks=false,breaklinks=true]{hyperref}
\usepackage{float}
\usepackage{graphicx}
\usepackage{color}

\usepackage{epigraph}
\definecolor{myurlcolor}{rgb}{0,0,0.7}
\usepackage{hyperref}
\hypersetup{colorlinks,
linkcolor=myurlcolor,
citecolor=myurlcolor,
urlcolor=myurlcolor}
\usepackage[mathscr]{euscript}

\input xy
\xyoption{all}

\newcommand{\maps}{\colon}    
\renewcommand{\H}{\mathcal{H}}  

\newcommand{\GL}{\mathrm{GL}}
\newcommand{\SL}{\mathrm{SL}}

\newcommand{\Z}{{\mathbb Z}}  
\newcommand{\C}{{\mathbb C}}  

\newcommand{\define}[1]{{\bf \boldmath{#1}}}

\theoremstyle{definition}

        \newcommand{\be}{\begin{equation}}
        \newcommand{\ee}{\end{equation}}
        \newcommand{\ba}{\begin{eqnarray}}
        \newcommand{\ea}{\end{eqnarray}}
        \newcommand{\ban}{\begin{eqnarray*}}
        \newcommand{\ean}{\end{eqnarray*}}
        \newcommand{\barr}{\begin{array}}
        \newcommand{\earr}{\end{array}}

\begin{document}
\title{The Moduli Space of Acute Triangles}
\author[Baez]{John C.\ Baez} 
\address{Department of Mathematics, University of California, Riverside CA, 92521, USA}
\date{October 1, 2023}
\maketitle

\begin{tikzpicture}
\node [] at (0, 0) {\includegraphics[width = 40 em]{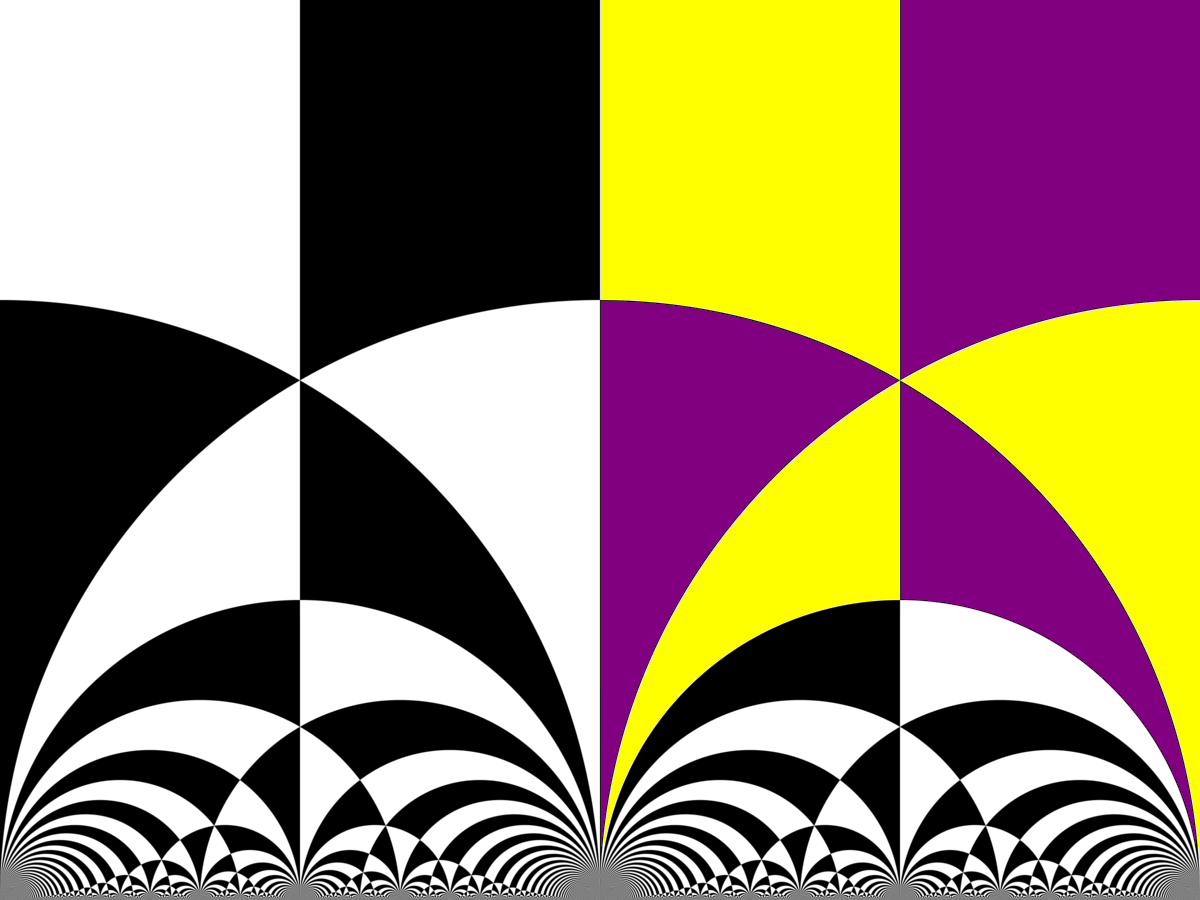}};
\node at (-7,-5.5) {$-1$};
\node at (-3.5,-5.5) {$-1/2\;$};
\node at (0,-5.5) {$0$};
\node at (3.5,-5.5) {$1/2$};
\node at (7,-5.5) {$1$};
\end{tikzpicture}
\vskip 1em

In mathematics we often like to classify objects up to isomorphism.  Sometimes the classification is discrete, but sometimes we have a notion of when two objects are `close'.  Then we can make the set of isomorphism classes into a topological space called a `moduli space'.   A simple example involves acute triangles.   We can define two triangles to be isomorphic if they are similar.  Then the moduli space of acute triangles is just the space of similarity classes of such triangles.

As a step toward explicitly describing this moduli space, first consider triangles in the complex plane with labeled vertices---that is, with a specified `first', `second' and `third' vertex.   Every such triangle is similar to one with the first vertex at $0$, the second at $1$, and the third at some point in the upper half-plane.  This triangle is acute precisely when its third vertex lies in this set:
\[      T = \left\{z \in \C \; \left\vert \; \mathrm{Im}(z) > 0, \; 0 < \mathrm{Re}(z) < 1, \; |z - \tfrac{1}{2}| > \tfrac{1}{2} \right.\right\}   \]
So, we say $T$ is the moduli space of acute triangles with labeled vertices.  This set
is colored yellow and purple above; the yellow and purple regions on top extend infinitely upward.

To get the moduli space of acute triangles with unlabeled vertices, we must mod out $T$ by the action of $S_3$ that permutes the three vertices.  To act by a permutation we relabel the vertices of the corresponding triangle, apply a similarity so that the new first and second vertices lie at 0 and 1, and then record the new location of the third vertex.  The six yellow and purple regions in $T$ are `fundamental domains' for this $S_3$ action: that is, they each contain exactly one point from each orbit.   If we reflect a labeled triangle corresponding to a point in a yellow region across the line $\mathrm{Re}(z) = \frac{1}{2}$, we get a triangle corresponding to a  point in a purple region, and vice versa.  Points on the boundary between two regions correspond to isosceles triangles.  All six regions meet at the point that corresponds to an equilateral triangle. 

The moduli space of acute triangles is closely related to a more famous moduli space: that of elliptic curves.  An \define{elliptic curve} is a torus equipped with the structure of a complex manifold.  But for those unfamiliar with complex manifolds, we can describe elliptic curves more concretely as follows.  Given a parallelogram in the complex plane we can identify its edges and get a torus called an elliptic curve.  Two parallelograms are said to give isomorphic elliptic curves if we can get one parallelogram by translating, rotating and/or dilating the other.  Note that we do not include reflections here: for example, applying complex conjugation to some parallelogram can yield a parallelogram giving a non-isomorphic elliptic curve.

For any point $z \in \H$ we can form a parallelogram with vertices $0, 1, z$ and $z+1$.    If we identify the opposite edges of this parallelogram we get an elliptic curve.  We can get every elliptic curve this way, up to isomorphism.  So, it is interesting to ask when two points in $\H$ give isomorphic elliptic curves.   

To answer this question, we introduce the group $\GL(2,\Z)$, consisting of invertible $2 \times 2$ integer matrices, and its subgroup $\SL(2,\Z)$ consisting of those matrices with determinant 1.   Both these groups act on the upper half-plane
\[     \H = \{z \in \C \; \vert \; \mathrm{Im}(z) > 0 \} \]
as follows: 
\[    \left( \begin{array}{cc} a & b \\ c & d \end{array} \right) \maps z \mapsto \frac{a z + b}{c z + d} .\]
Each of the light or dark regions shown above is a fundamental domain for the action of $\GL(2,\Z)$.   Elements of $\GL(2,\Z)$ with determinant $-1$ map light regions to dark ones and vice versa.   Elements with determinant $1$ map light regions to light ones and dark ones to dark ones.   People more often study the action of the subgroup $\SL(2,\Z)$.  The union of any light region and dark one sharing an edge forms a fundamental domain of $\SL(2,\Z)$.  

We can now answer our question: two points $z,z' \in \H$ give isomorphic elliptic curves if and only if  $z' = g z$ for some $g \in \SL(2,\Z)$.   Thus the quotient space $\H/\SL(2,\Z)$ is the \define{moduli space of elliptic curves}: points in this space correspond to isomorphism classes of elliptic curves.   For details, see any good introduction to elliptic curves \cite{Knapp,Koblitz}.

Since $T$ is the union of six fundamental domains for $\GL(2,\Z)$---the yellow and purple regions in the figure---it is the union of three fundamental domains for $\SL(2,\Z)$.  There is thus a map
\[          p \maps T \to \H/\SL(2,\Z) \]
from the moduli space of acute triangles with labeled vertices to the moduli space of elliptic curves, and generically this map is three-to-one.  This map is not onto, but if we take the closure of $T$ inside $\H$ we get a larger set
\[      \overline{T} = \left\{z \in \C \; \left\vert \; \mathrm{Im}(z) > 0, \; 0 \le \mathrm{Re}(z) \le 1, \; |z - \tfrac{1}{2}| \ge \tfrac{1}{2} \right.\right\}  \]
whose boundary consists of points corresponding to right triangles.  Then $p$ 
extends to an onto map
\[       p \maps \overline{T} \to \H/\SL(2,\Z) .\]

The existence of this map suggests that from any acute or right triangle in the plane we can construct an elliptic curve.   This is in fact true.   How can we understand this more directly?  

Take any acute or right triangle in the complex plane.  Rotating it 180${}^\circ$ around the midpoint of any edge we get another triangle.  The union of these two triangles is a parallelogram.   Identifying opposite edges of this parallelogram, we get an elliptic curve! There are three choices of how to build this parallelogram, one for each edge of the original triangle, but they give isomorphic elliptic curves.   Moreover, two acute or right triangles differing by a translation, rotation or dilation give isomorphic elliptic curves.  Even better, every elliptic curve is isomorphic to one arising from this construction.  Thus, this construction gives a map from $\overline{T}$ onto $\H/\SL(2,\Z)$, and with a little thought one can see that this map is $p$.

I learned about the moduli space of acute triangles from James Dolan.  There has also been interesting work on the moduli space of \emph{all} triangles in the plane.  Gaspar and Neto  \cite{GasparNeto} noticed that this space is a triangle, and Stewart later gave a more geometrical explanation \cite{Stewart}.  In fact all the moduli spaces mentioned here are better thought of as moduli `stacks': stacks give a way to understand the special role of more symmetrical objects, like isosceles and equilateral triangles.  Behrend \cite{Behrend} has written an introduction to stacks using various moduli stacks of triangles and the moduli space of elliptic curves as examples.   Though he does not describe the map $p$ or its stacky analogue, his work is a nice way to dig deeper into some of the ideas discussed here.

\subsubsection*{Acknowledgements} 
I thank Roice Nelson for the figure here, and James Dolan for many conversations on this subject.

\end{document}